\documentstyle[12pt]{article}
\def\eqnum#1{}

\begin{document}

\bigskip
\centerline{\bf ON THE UNIQUENESS OF THE BRANCHING}
\centerline{\bf  PARAMETER FOR  A RANDOM CASCADE MEASURE}

\bigskip
\centerline{\bf G. Molchan$^{1,2}$}
\centerline{$^1$ \it Observatoire de la Cote d'Azur}
\centerline{\it B.P 4229, 06304, Nice Cedex 4, France}
\centerline{$^2$ \it International Institute of Earthquake Prediction Theory}
\centerline{\it Russian Academy of Sciences,}
\centerline{\it Warshavskoye sh.79, k.2, Moscow 117556, Russian Federation,}
\centerline{\it e-mail: molchan@mitp.ru}
\qquad

\begin{abstract}

\noindent {\small
An independent random cascade measure $\mu$ is specified by a random
generator $(w_1,...,w_c)$, $E\sum w_i=1$
where $c$ is the branching parameter. It is shown under certain
restrictions that,
if  $\mu$  has two generators with a.s. positive components, and the
ratio $\ln\,c_1/\ln\,c_2$ for their branching parameters is an irrational
number, then $\mu$ is a Lebesgue measure. In other
words, when $c$ is a power of an integer number $p$ and the $p$ is minimal
for $c$, then a cascade measure that has the property of intermittency
specifies $p$ uniquely.} \\

\end{abstract}

\newpage

\bigskip

{\bf INTRODUCTION}

\bigskip
Intermittency in turbulence is usually expressed in scaling terms of
structure functions. For energy $\varepsilon (\Delta_x)$ dissipated in
a cell $\Delta_x$, empirical data give

\begin{eqnarray}
<\varepsilon (\Delta_x)^q>\,\,  \sim \,\, |\Delta|^{\tau(q)+d},
\quad |\Delta|\ll 1
\eqnum{1}
\end{eqnarray}
where $<\cdot >$ denotes spatial averaging, $|\Delta|$ is the linear
cell size, and $d$ is the spatial dimension ($d=1$ in what follows).
Intermittency corresponds to scaling exponents $\tau (q)$ of $a$ nonlinear
type for $q>0$. Historically the first interpretation of
intermittency is associated with
Richardson's idea (see, e.g., ref. 1) as to energy being transmitted
from larger to smaller scales in an inertial range $(L, \delta)$: $L$ is the
external scale, while $\delta \ll L$ is the Kolmogorov scale at which
the dissipation occurs.

The above idea can be formalized by means of the following
recursive procedure which defines an independent random cascade.
We denote by $\varepsilon^{(n)}_\alpha$ the energy in a cell
$\Delta^{(n)}_\alpha$ of level $n$.
Each cell $\Delta_\alpha^{(n)}$ is divided into $c$ equal subcells of
level $(n+1)$;
into these the energy $\varepsilon_\alpha^{(n)}$ is transmitted with
random coefficients

\begin{eqnarray}
(w_1(\Delta_\alpha^{(n)}),...,w_c(\Delta_\alpha^{(n)}):
= W(\Delta_\alpha^{(n)}).
\nonumber
\end{eqnarray}

The vectors $W(\Delta_\alpha^{(n)})$, which are called breakdown
coefficients in the physics literature, are statistically independent
and identically distributed for different cells
$\Delta_\alpha^{(n)}$ of all levels $n$. Their distribution is
specified by the random vector (or {\it cascade generator})
$W=(w_1,...,w_c)$, for which $w_i\ge 0$ and $E\sum w_i=1$ , corresponding
to the law of conservation of
energy in the average.
In what follows  we will restrict our consideration to cascade
generators for which $P(w_*=1)<1$ and $P(w_*>0)=1$. Here $w_*$ is
the normalized random component of the vector $W$, namely

\begin{eqnarray}
w_* = \{ cw_i \quad \mbox {\rm with probability} \quad 1/c
\nonumber
\end{eqnarray}
and $Ew_*=1$.

On can express many cascade properties in terms of $w_*$.
In particular, the following condition: $Ew_*\log_cw_*<1$
ensures the existence of a nontrivial limit of
the measures

\begin{eqnarray}
\varepsilon_n(dx) = \sum_\alpha  \varepsilon_\alpha^{(n)}
{\bf 1}_{\Delta_\alpha^{(n)}}(x)\,dx/|\Delta^{(n)}|
\nonumber
\end{eqnarray}
as $n\to \infty$ (refs. 2,3). Following Mandelbrot (ref. 4),
the limiting cascade measure $\varepsilon (dx)$ is considered as
a  model of dissipated energy field in turbulence.

Under very general conditions the cascade measure
$\varepsilon (dx)$ has the intermittency property (1)
where  $\sim$  denotes logarithmic asymptotics, i.e., $a\sim b$
when $\ln a=\log b(1+o(1))$ a.s. (refs. 5,6). The scaling
exponents in (1) are closely connected with the
function

\begin{eqnarray}
\tau^H(q) = q - \log_cEw_*^q - 1,
\eqnum{2}
\end{eqnarray}
which represents the heuristic estimate of $\tau (q)$. It is easily
found by replacing $\varepsilon (dx)$ with the pre-limit measure
$\varepsilon_n(dx)$, $n\gg 1$ and the spatial averaging $<\cdot >$ with
ensemble averaging, i.e., with the operation of mathematical
expectation $E$. With large $|q|$ these manipulations lead to false
estimates of $\tau$ (ref. 5). The true function $\tau (q)$ is identical
with (2) in the interval $q\in (q_-, q_+)$, only where
$-\infty \le  q_- \le 0$ and $1 \le q_+ \le \infty$.
The function $\tau$ is linear outside of $(q_-,q_+)$
if $q_-<0<1<q_+$: $\tau (q)=a_{\pm}q$. Both lines
$a_{\pm}q$ for finite $q_{\pm}$ are tangent to the curve of $\tau^H(q)$,
which uniquely specifies the critical points  $q_-$ and $q_+$ as
those tangent points closest to 0.

The interpretation of intermittency in terms of cascades uses two
assumptions that are not particularly attractive from the standpoint
of physics:

(a) the ratio of adjacent scales $|\Delta_n|/|\Delta_{n+1}|$, i.e.,
the cascade's branching parameter, is fixed;

(b) the Kolmogorov dissipation scale is zero: $\delta=0$.

>From (a) it follows that $L/\delta =c^N$, where $c$ and $N$ are integers.
Varying $L$, which is natural for many physical objects, we vary $c$ thereby.
For this reason it is desirable to deal with cascades whose
statistical properties are independent of the parameter $c$.
This standpoint has proved fruitful for resolving the problem
of parametrization of empirical $\tau$-functions. A broad class of
functions (2) corresponding infinitely divisible random
variables  $\log w_*$ was suggested
for practical purposes
(refs. 7-11). Any $\tau$-function of this type
can be produced by a cascade
generator of any dimension. Unfortunately, a complete description
of $\tau$-functions that would
have the above property is unknown.

In some applications there are  attempts to introduce a
scale densification in the cascade process [refs. 10, 12].
The aim of this modification is twofold: to get rid of the
above assumptions (a,b) and to justify
a "universal class" of cascades. This idea has unfortunately
remained without justification.

An opposite viewpoint on the parameter $c$
for turbulent cascades is contained
in ref. 13 where the authors assumed $c=2$, since  the Navier-Stokes
equation involved a nonlinearity of the second order. Based on this
assumption, the authors derive the statistical conclusion that the
coefficients $W(\Delta^{(n)})$ are interdependent for two adjacent
levels $n$ and $(n+1)$ in actual turbulence. The conclusion lacks
experimental corroboration of the assumed hypothesis $c=2$. Otherwise
it can equally well be regarded as an artefact.

It is our purpose to show that, under conditions that are
natural for turbulence, the least integer-valued parameter $p$ in the
representation $c=p^n$, $n\ge 1$ is uniquely specified by a cascade measure
($c=p=2$ in the case of ref. 13). From this it follows that a locally
positive
cascade measure having the intermittency property and a two
generators of
significantly different dimensions, i.e., when
$\log c_1/\log c_2$ is irrational, does not exist.
In other words, the requirement that the cascade measure be independent
of the branching parameter is much too fine for the phenomenological
model of intermittency. However, if the cascade measure is regarded
as the model of a physical object, the above parameter $p$ in $c=p^n$ should
have a physical meaning, hence an algorithm is required to identify
it from the cascade measure. Such
algorithms are unknown to us.

The present study generalizes the results of my previous work
(ref. 14).

\bigskip

{\bf THE MAIN RESULT}

\bigskip
This section consists in the following

{\bf Theorem.} Suppose a random cascade measure $\mu$ on $[0,1]=I$ is
locally positive, i.e., $\mu(\Delta)>0$ a.s. for any subinterval
$\Delta \subset I$, the total mass $M=\mu(I)$ has a second moment,
$EM^2<\infty$, and  $q_+>2$. If $\mu$ has
two generators $\xi \in R^{c_1}$ and $\eta \in R^{c_2}$,
$0<c_1<c_2$, and $\log c_1/\log c_2$ is irrational, then $\mu$ is
a Lebesgue measure.

Let us comment on the conditions of this theorem. The main
requirements, namely, that $\mu$ should be locally positive and
$\log c_1/\log c_2$ should be irrational are essential.
For instance, the cascade generator $W$ and the tensor product of
its independent copies $W_1\otimes W_2$ generate
the same cascade measure having the branching parameters $c$ and $c^2$.
A measure of the type $\mu (dx)=\delta (x-\xi)dx$ where $\xi$ is a
random uniformly distributed variable on [0,1] is a
cascade measure with a generator of arbitrary dimension:
$W=(0,...,0,1,0,...0)$. Here, 1 occupies
the $i$-th position with probability $1/c$.

The requirements $EM^2<\infty$ and $q_+>2$ are purely technical in character
and are merely needed in the method of proof we employ.
Under these conditions $\tau(q)=\tau^H(q)$ for $0<q<2$.
We remind that
the tangent to $\tau^H(q)$ at the point $q_+<\infty$ passes through (0,0).
To be more specific, the function $-\tau^H(q)$ is convex, so that one
should speak of the support line at the point $q_+<\infty$ rather than of
the tangent. Judging by empirical evidence (see, e.g., ref. 1),
these requirements do not constitute restrictions on turbulent cascades.

The proof rests on two Lemmas.

{\bf Lemma 1.} Suppose two vectors  $\xi =(\xi_0,...,\xi_{c_1-1})$ and
$\eta =(\eta_0,...,\eta_{c_2-1})$, $1<c_1<c_2$ with
positive components commute with respect to the tensor
product: $\xi \otimes \eta =\eta \otimes \xi$.
If $\log c_1/\log c_2$ is irrational, then both vectors have
constant components.

{\it Proof of Lemma 1}. We write down the commutation condition
for the vectors $\xi$ and $\eta$ as follows;

\begin{eqnarray}
\eta_{[q]_{c_2}} \xi_{\{ q\}_{c_2}} = \xi_{[q]_{c_1}}\eta_{\{ q\}_{c_1}},
\quad 0 \le q < c_1c_2,
\eqnum{3}
\end{eqnarray}
where $[q]_n$ and $\{ q\}_n$ are the integer part and the remainder
resulting from dividing $q$ by $n$. One has for $q=\alpha <c_1$:

\begin{eqnarray}
\eta_0\xi_\alpha = \xi_0\eta_\alpha, \quad 0 \le \alpha < c_1.
\nonumber
\end{eqnarray}

Consequently, one can assume $\xi_\alpha =\eta_\alpha$,  $\alpha <c_1$
without loss of generality; (3) then becomes

\begin{eqnarray}
\eta_{[q]_{c_2}} \eta_{\{ q\}_{c_2}} = \eta_{[q]_{c_1}}\eta_{\{ q\}_{c_1}},
\quad 0 \le q < c_1c_2.
\eqnum{4}
\end{eqnarray}

In particular,

\begin{eqnarray}
\eta_{pc_1+\alpha} = \eta_p\eta_\alpha /\eta_0,
\quad 0 \le pc_1+\alpha < c_2, \quad 0 \le  \alpha <c_1.
\eqnum{5}
\end{eqnarray}

Iteration yields

\begin{eqnarray}
\eta_\beta/ \eta_0 = \prod_{i=0}^k(\eta_{\alpha_i} /\eta_0),
\quad \beta = \alpha_0+\alpha_1c_1+...+\alpha_kc_1^k < c_2,
\quad 0 \le  \alpha_i < c_1,
\eqnum{6}
\end{eqnarray}
that is, the vector $\eta$ can be uniquely reconstructed from the
first $c_1$ coordinates.

Suppose $D$ is the greatest common divisor of $c_1$, $c_2$ and
$D:=(c_1, c_2)<c_1$. One can then find integer
$\alpha_0<c_1/D:=k_1$ and $\beta_0<c_2/D:=k_2$ such that
$\alpha_0c_2=\beta_0c_1+D$.

Denote

\begin{eqnarray}
\nonumber
a_p&=&\eta_{\beta_0+p}/ \eta_{\alpha_0}, \quad \beta_0+p<c_2,\\
\nonumber
\vec{\eta}_r&=&(\eta_{rD}, \eta_{rD+1},...,\eta_{rD+D-1}),
\quad 0 \le r < k_2.
\nonumber
\end{eqnarray}

Since $\alpha_0<c_1$, one has
$\alpha_0c_2+\beta =\beta_0c_1+D+\beta <c_1c_2$  for any
$0\le \beta <c_2$. Hence, using (3) and
the above notation, one has:

\begin{eqnarray}
\eta_{\,_\beta} = a_{\,_{[\beta + D]_{c_1}}}  \eta_{\,_{\{ \beta +D\}_{c_1}}},
\quad 0 \le \beta <c_2.
\eqnum{7}
\end{eqnarray}

>From (7) it follows that

\begin{eqnarray}
\vec{\eta}_{pk_1+i} = \left\{\begin{array}{ll}
a_p\vec{\eta}_{i+1},\quad 0 \le  i < k_1-1 \\
a_{p+1}\vec{\eta}_0,\quad i = k_1-1,
\end{array}\right.
\eqnum{8}
\end{eqnarray}
where $pk_1+i<k_2$. Put $p=0$ here. One then arrives at the
recurrence relation

\begin{eqnarray}
\vec{\eta}_i = a_0\vec{\eta}_{i+1}, \quad 0 \le i < k_1-1,
\eqnum{9}
\end{eqnarray}
whence

\begin{eqnarray}
\vec{\eta}_i = a_0^{-i}\vec{\eta}_0, \quad 0 \le i < k_1.
\eqnum{10}
\end{eqnarray}

>From (5) one  has

\begin{eqnarray}
\vec{\eta}_{pk_1+i} = (\eta_p/\eta_0)\vec{\eta}_i, \quad 0 \le i < k_1,
\quad pk_1+i < k_2.
\eqnum{11}
\end{eqnarray}

The use of (8, 11, 9) yields the chain of relations

\begin{eqnarray}
a_p\vec{\eta}_{i+1} \stackrel{\rm (8)}{=}
\vec{\eta}_{pk_1+i} \stackrel{\rm (11)}{=}
(\eta_p/\eta_0)\vec{\eta}_i \stackrel{\rm (9)}{=}
a_0(\eta_p/\eta_0)\vec{\eta}_{i+1}
\quad 0 \le i < k_1-1,
\nonumber
\end{eqnarray}
where $0\le pk_1+i<k_2$. Put $i=0$ here. Then

\begin{eqnarray}
a_p/a_0 = \eta_p/\eta_0, \quad 0 \le pk_1 < k_2.
\eqnum{12}
\end{eqnarray}

We now make use of (8) with $i=k_1-1$:

\begin{eqnarray}
\nonumber
a_{p+1}\,\vec{\eta}_0 \stackrel{\rm (8)}{=}
\vec{\eta}_{pk_1+k_1-1} \stackrel{\rm (11)}{=}
(\eta_p/\eta_0)\,\vec{\eta}_{k_1-1} \stackrel{\rm (10)}{=}
(\eta_p/\eta_0)a_0^{-k_1+1}\,\vec{\eta}_0, \\
\nonumber
\quad 0\le pk_1+k_1-1<k_2,\qquad \qquad \qquad \qquad
\nonumber
\end{eqnarray}
i.e.,

\begin{eqnarray}
a_{p+1}/a_0 = (\eta_p/\eta_0)a_0^{-k_1} \stackrel{\rm (12)}{=}
(a_p/a_0)a_0^{-k_1}, \quad 0\le pk_1+k_1-1<k_2,
\nonumber
\end{eqnarray}
whence

\begin{eqnarray}
a_p/a_0 = a_0^{-pk_1}, \quad 0 \le p<(k_2+1)/k_1,
\eqnum{13}
\end{eqnarray}

\begin{eqnarray}
\eta_p/\eta_0 = a_0^{-pk_1}, \quad 0 \le p<k_2/k_1.
\eqnum{14}
\end{eqnarray}

>From (10, 11, 14) one has for $r=pk_1+i<k_2$, $i<k_1$:

\begin{eqnarray}
\vec{\eta}_r \stackrel{\rm (11)}{=}
(\eta_p/\eta_0)\,\vec{\eta}_i \stackrel{\rm (10)}{=}
(\eta_p/\eta_0)\,a_0^{-i}\,\vec{\eta}_0 \stackrel{\rm (14)}{=}
a_0^{-r}\vec{\eta}_0,
\nonumber
\end{eqnarray}
i.e.,

\begin{eqnarray}
\vec{\eta}_r = a_0^{-r}\vec{\eta}_0, \quad 0\le r<k_2.
\eqnum{15}
\end{eqnarray}

The original relations (4) when expressed in terms of the $\vec{\eta}_i$
have the form

\begin{eqnarray}
\eta_{[r]_{k_1}}\vec{\eta}_{\{ r\}_{k_1}} =
\eta_{[r]_{k_2}}\vec{\eta}_{\{ r\}_{k_2}},
\quad 0\le r<c_1c_2/D.
\nonumber
\end{eqnarray}
Hence in virtue of (15) one has

\begin{eqnarray}
\eta_{[r]_{k_1}}a_0^{-\{ r\}_{k_1}} =
\eta_{[r]_{k_2}}a_0^{-\{ r\}_{k_2}}.
\eqnum{16}
\end{eqnarray}

Let $r=k_2=pk_1+i$, $i<k_1$. Then $p<k_2/k_1$ and

\begin{eqnarray}
\eta_pa_0^{-i} = \eta_1
\eqnum{17}
\end{eqnarray}
or
\begin{eqnarray}
a_0^{-pk_1} \stackrel{\rm (14)}{=}
\eta_p/\eta_0 \stackrel{\rm (17)}{=}
(\eta_1/\eta_0)a_0^i \stackrel{\rm (16)}{=}
a_0^{-k_1+i}.
\nonumber
\end{eqnarray}
Hence $a_0=1$. Otherwise $pk_1=k_1-i$ or $k_2=pk_1+i=k_1$,
which is impossible.

>From (16) and $a_0=1$ one gets

\begin{eqnarray}
\eta_{[r]_{k_1}} = \eta_{[r]_{k_2}}.
\nonumber
\end{eqnarray}

One has

\begin{eqnarray}
\eta_p=\eta_0, \quad 0\le p<k_2/k_1
\nonumber
\end{eqnarray}
when $r=pk_1+i<k_2$ and

\begin{eqnarray}
\eta_p=\eta_1, \quad k_2/k_1 < p < 2k_2/k_1
\nonumber
\end{eqnarray}
when $k_2<pk_1+i<2k_2$
and so $\eta_p=\eta_0$, $p<2k_2/k_1$, since $\eta_1=\eta_0$.

Proceeding as above, we shall prove in succession that
the components of $\vec{\eta}_0$ are constant.
By (15) $\vec{\eta}_r = \vec{\eta}_0$,  $0\le r<k_2$, thus
we have the desired relation
$\eta =(\eta_0,...,\eta_0)$.

Consider the case in which $c_2$ is divisible by $c_1$,
i.e., $D=c_1$. Put $q=rc_1$, $r<c_2$ in (4). One gets

\begin{eqnarray}
\eta_{[r]_{k_2}}\eta_{c_1\{ r\}_{k_2}} = \eta_r\eta_0
\stackrel{\rm (5)}{=}
\eta_{[r]_{c_1}}\eta_{\{ r\}_{c_1}},
\eqnum{18}
\end{eqnarray}
where $k_2=c_2/c_1$. However, if $c_2$ is divisible by $c_1$, then by (4),

\begin{eqnarray}
\eta_p = \eta_{pc_1}, \quad 0\le p<k_2,
\nonumber
\end{eqnarray}
i.e., $\eta_{c_1\{ r\}_{k_2}}=\eta_{\{ r\}_{k_2}}$.
Consequently, (18) means that the vectors \quad
$\xi = (\eta_0,\eta_1,...,\eta_{c_1-1})$
and $\eta =(\eta_0,...,\eta_{k_2-1})$  commute with respect to the
tensor product. The maximum dimension of the new vectors
$\xi$ and $\eta$ has decreased from $c_2$ to
$\max (c_1,k_2)$. The problem has reduced to the case
already considered. The process of reducing the dimension of $\xi ,\eta$ is
finite, terminating when the dimensions $c_1^{(k)}$ and $c_2^{(k)}$ are not
mutually divisible at the stage $k$. Otherwise, as is easily seen,
$c_1=p^{k_1}$ and $c_2=p^{k_2}$, where $p,k_1,k_2$ are integers. That
case is ruled out, because  $\ln c_1/\ln c_2$
is irrational. When $c_1^{(k)}$ and $c_2^{(k)}$ are not mutually divisible,
induction on $k$ applied to (15) will demonstrate that $\eta$ has constant
components. Lemma 1 is proven. $\diamond$

Lemma 1 yields an immediate corollary which will be stated
here as Lemma 2.

{\bf Lemma 2.} Let a random cascade measure $\mu$ on $I=[0,1]$ has
two random generators $\xi =(\xi_0,...,\xi_{c_1-1})$ and
$\eta =(\eta_0,...,\eta_{c_2-1})$,  $0<c_1<c_2$, and $\ln c_1/\ln c_2$ be
irrational. If the measure $\mu$ is locally positive a.s. and
$E[\mu (I)]^\rho <\infty$  for some $\rho >1$, then
$E\xi_\alpha^\rho = E\xi_0^\rho$, $0\le \alpha <c_1$ and
$E\eta_\beta^\rho = E\eta_0^\rho$, $0\le \beta <c_2$. Also,

\begin{eqnarray}
\biggr [\frac{E\xi_0^{\rho}}{(E\xi_0)^\rho} \biggl ]^{\frac{1}{\ln c_1}} =
\biggr [\frac{E\eta_0^{\rho}}{(E\eta_0)^\rho} \biggl ]^{\frac{1}{\ln c_1}},
\eqnum{19}
\end{eqnarray}
provided $0<\rho <q_+$.

{\it Proof of Lemma} 2.

It follows from the definition of the random cascade
measure $\mu$ (ref. 4) that it satisfies the following stochastic
equation:

\begin{eqnarray}
\mu \stackrel{\rm (d)}{=} \sum z_i\mu_i\circ T_i^{-1}
\eqnum{20}
\end{eqnarray}
where $z=(z_0,...,z_{c-1})$ is the generator of $\mu$, the $\mu_i$ are
independent copies of $\mu$ that are statistically independent of $z$
as well, and $T_ix=(i+x)/c$ is a linear mapping of the interval $I$
into $I_i=[i/c, (i+1)/c)$. Note that
$\mu \circ T_i^{-1}(I_j)=0$ for $i\ne j$. The equality
$\mu_1 \stackrel{\rm (d)}{=} \mu_2$ for the two measures
means that the distributions of $\mu_i(f)$, $i=1,2$ are equal
for any smooth finite functions $f:R^1\to R^1$.

Suppose $\xi$ and $\eta$ are the generators of $\mu$.
Use (20) with weights $z=\xi$, and then use the same representation
for each measure $\mu_i$ with weights $\eta^{(i)}$.
The weights $\eta^{(i)}$ are independent copies of $\eta$ which are
independent of $\xi$. The result is

\begin{eqnarray}
\mu \stackrel{\rm (d)}{=} \sum_{0\le \alpha <c_1} \xi_\alpha
\sum_{0\le \beta <c_2} \eta^{(\alpha)}_\beta
\mu_{\alpha,\beta}\circ T_{\alpha,\beta}^{-1}\,,
\eqnum{21}
\end{eqnarray}
where $\mu_{\alpha,\beta}$ are independent copies of $\mu$ that are also
independent of $\xi$, $\eta^{(\alpha)}$, $\alpha <c_1$, while
$T_{\alpha,\beta}$ is a linear mapping of $I$ into the interval
$\delta_{\alpha,\beta}=[c_2\alpha +\beta,\, \, c_2\alpha +\beta +1]/(c_1c_2)$,
$0\le \alpha <c_1$, $0\le \beta <c_2$.
Interchanging $\xi$ and $\eta$ in (21), one gets
a representation of $\mu$ that involves the random quantities
$\{ \tilde{\eta}_{\beta'}\}$ and $\{ \tilde{\xi}_{\alpha'}^{(\beta')}\}$
and the intervals
$\delta_{\beta',\alpha'}=[c_1\beta' +\alpha',\,\,
c_1\beta' +\alpha' +1]/(c_1c_2)$.
Obviously, one has $\delta_{\alpha,\beta}=\delta_{\beta',\alpha'}$
when $c_2\alpha +\beta =c_1\beta'+\alpha'$. Consequently,

\begin{eqnarray}
E[\mu (\delta_{\alpha,\beta})]^\rho
&=& E[\xi_\alpha \eta^{(\alpha)}_\beta
\mu_{\alpha,\beta}(\delta_{\alpha,\beta})]^\rho =
E\xi_\alpha^\rho E\eta_\beta^\rho m_\rho \\
\nonumber
&=& E[\eta_{\beta'}\xi_{\alpha'}^{(\beta')}
\mu_{\alpha',\beta'}(\delta_{\beta',\alpha'})]^\rho =
E\xi_{\alpha'}^\rho E\eta_{\beta'}^\rho m_\rho,
\eqnum{22}
\end{eqnarray}
where $m_\rho =E\mu^\rho(I)$. According to (ref. 15), the requirement
$0<m_\rho <\infty$, $\rho >1$ is equivalent to
$E\xi_\alpha^\rho <\infty$ ($E\eta_\beta^\rho <\infty)$
for all components of the generator. It also follows from (22)
that the moments are $E\xi_\alpha^\rho >0$, because
$\mu (\Delta )>0$ a.s. for any subinterval of $I$.

Relation (22) means that the vectors
$\{ E\xi_\alpha^\rho, \alpha =0,...,c_1-1\}$ and
$\{ E\eta_\beta^\rho,\, \beta =0,...,c_2-1\}$ are
positive and commute with respect to the tensor product.
The use of Lemma 1 therefore yields the right-hand side
of Lemma 2:
$E\xi_\alpha^\rho = E\xi_0^\rho$, $0\le \alpha <c_1$,
$E\eta_\beta^\rho = E\eta_0^\rho$,  $0\le \beta <c_2$.
If $0<\rho <q_+$, one can equate the $\tau^H(\rho)$ for
the generators  $\xi$ and $\eta$. From (2) one has

\begin{eqnarray}
\log_{c_1}E\xi_*^\rho = \log_{c_2}E\eta_*^\rho.
\eqnum{23}
\end{eqnarray}
But $E\xi_*^\rho = E\sum\limits_{\alpha} \xi_\alpha^\rho c_1^{\rho -1}$
with  $E\xi_\alpha^\rho = E\xi_0^\rho$  and  $E\xi_\alpha = E\xi_0=c^{-1}$.
For this reason one has $E\xi_*^\rho = E\xi_0^\rho /(E\xi_0)^\rho$.
Similarly, $E\eta_*^\rho = E\eta_0^\rho /(E\eta_0)^\rho$.
Substitution of these relations in (23) yields (19).
$\diamond$

The proof of the theorem uses another obvious number-theoretic
fact which will be treated as a separate statement.

{\bf Statement 3.} Suppose that the integer numbers $n_1$
and $n_2$ are not mutually divisible, and that
$T\alpha =\{ n_2\alpha \}_{n_1}$
where $\{ k\}_n$ is the remainder left after dividing  $k$ by $n$.
One can then find $0<\alpha <n_1$ and $k(\alpha)>0$ such that
$T^{k(\alpha)}\alpha =\alpha$.

{\bf Proof of the Theorem}.

We are going to make use of two representations of $\mu$ in the form
(21). The one is based on the independent generators
$\xi, \eta^{(\alpha)},  \alpha =0,...,c_1-1$
and the other on $\eta$ and $\xi^{(\beta)}, \beta =0,...c_2-1$.
Consider the values of $\mu$ on elements of the partitioning
${\cal F}$ of [0,1] into $c_1c_2$ equal parts.
One then gets equality of distribution for two families of

\begin{eqnarray}
\{ \xi_\alpha \eta_\beta^{(\alpha)}M_{\alpha \beta}\}
\stackrel{\rm (d)}{=}
\{ \eta_{\beta'}\xi_{\alpha'}^{(\beta')}M_{\alpha' \beta'}\},
\nonumber
\end{eqnarray}
where $0\le \alpha,\, \alpha' <c_1$, $0\le \beta,\, \beta' <c_2$ and
$M_{\alpha \beta}$ are independent copies of $\mu (I)$ under
the following correspondence between the subscripts:

\begin{eqnarray}
q(\alpha, \beta):= c_2\alpha + \beta = c_1\beta' + \alpha'.
\nonumber
\end{eqnarray}
Here, $q$, $0\le q<c_1c_2$  is the natural numbering of elements of
the partioning ${\cal F}$ on [0,1]. As follows from Lemma 2, the moments

\begin{eqnarray}
E[\xi_\alpha \eta_\beta^{(\alpha)}M_{\alpha \beta}]^\rho =
E\xi_0^\rho \,E \eta_0^\rho m_\rho, \quad  \rho =1,2,
\nonumber
\end{eqnarray}
are independent of $\alpha, \beta$. We shall make use of that circumstance
for calculating the moments $E\mu (\delta_q)\mu (\delta_{q+1})$ where
$\delta_q$ is an element of ${\cal F}$ with index $q$. The twofold
representation (21) for $\mu$  yields a set of equations. We want
to write it down in compact form by first denoting

\begin{eqnarray}
\nonumber
a_\alpha &=& E \xi_{\alpha  -1}\xi_\alpha /E\xi_0^2,  \qquad
b_\beta  = E\eta_{\beta -1}\eta_\beta /(E\eta_0)^2 \\
\nonumber
\vec{a} &=& (a_1,...,a_{c_1-1}), \quad \, \,
\quad \vec{b} = (b_1,...,b_{c_2-1}) \\
\nonumber
V_a &=& E \xi_0^2 /(E\xi_0)^2,  \qquad \quad \,
V_b = E\eta_0^2 /(E\eta_0)^2
\nonumber
\end{eqnarray}
Now note that

\begin{eqnarray}
E\eta_{\beta_1}^{(\alpha_1)} \eta_{\beta_2}^{(\alpha_2)} =
\left\{\begin{array}{ll}
(E\eta_0)^2, \qquad \alpha_1 \ne \alpha_2 \\
E\eta_{\beta_1} \eta_{\beta_2},\quad \, \, \alpha_1 = \alpha_2.
\end{array}\right.
\nonumber
\end{eqnarray}

Similar equalities also hold for $\xi_\alpha^{(\beta)}$. Consequently,
the equations derived by calculating the moments
$E\mu (\delta_q) \mu (\delta_{q+1})$ in two different ways
have the form $X=Y$ where

\begin{eqnarray}
\nonumber
X = \{ V_b\vec{a};\,\,\,b_1,V_b\vec{a};\,\,\,b_2,V_b\vec{a};...;
b_{c_2-1},V_b\vec{a}\}, \\
\nonumber
Y = \{ V_a\vec{b};\,\,\,a_1,V_a\vec{b};\,\,\,a_2,V_a\vec{b};...;
a_{c_1-1},V_a\vec{b}\}.
\nonumber
\end{eqnarray}

It is our aim to show that the relation $X=Y$ yields
$V_a=E \xi_0^2 /(E\xi_0)^2=1$, i.e.,
the variance of $\xi_0$, hence that of $\xi_\alpha$,\,
$0\le \alpha <c_1$, equals zero.
Consequently, the generator $\xi$ has identical nonrandom
components, so that $\mu$ is a Lebesgue measure.

>From the fact that the first $c_1-1$ coordinates of $X$ and $Y$
are equal one has

\begin{eqnarray}
V_ba_i = V_ab_i, \quad 1\le i<c_1.
\nonumber
\end{eqnarray}

For this reason the coordinates of the vectors $X$ and $Y$ in the
notation $V_ab_i=\tilde{b}_1$ have the form

\begin{eqnarray}
X_r = \left\{\begin{array}{ll}
\,\,\tilde{b}_{\{ r\}_{c_1}}, \quad \quad  \{ r\}_{c_1} \ne 0
\qquad \qquad \mbox {\rm (a)} \\
V_a^{-1}\tilde{b}_{[r]_{c_1}}, \quad \{ r\}_{c_1}=0 \qquad
\qquad \mbox {\rm (b)}
\end{array}\right.
\eqnum{24}
\end{eqnarray}

\begin{eqnarray}
Y_r = \left\{\begin{array}{ll}
\,\,\tilde{b}_{\{ r\}_{c_2}}, \quad \quad  \{ r\}_{c_2} \ne 0
\qquad \qquad \mbox {\rm (a)} \\
V_b^{-1}\tilde{b}_{[r]_{c_2}}, \quad \{ r\}_{c_2}=0 \qquad
\qquad \mbox {\rm (b)}
\end{array}\right.
\eqnum{25}
\end{eqnarray}
Putting $r=c_1k$, one derives

\begin{eqnarray}
b_k = \left\{\begin{array}{ll}
Vb_{\{ c_1k\}_{c_2}}, \quad \quad \,\, \{ c_1k\}_{c_2} \ne 0 \\
Ub_{c_1k/c_2}, \quad \quad \, \, \, \, \{ c_1k\}_{c_2}=0
\end{array}\right.
\eqnum{26}
\end{eqnarray}
where $V=V_a$, $U=V_a/V_b$. Relation (26) also holds for the
$b_k$ and $\tilde{b}_k$.

Let $c_2=c_1^{r_1}c_2'$ where $r_1\ge 0$ is the maximum multiplicity
of $c_1$ in $c_2$.

{\it Consider the case} $(c_2',c_1)<\min (c_2',c_1)$.

Here, $c_2'>1$, because $\ln c_2/\ln c_1$ is
irrational. Put $k=c_1^{r_1}p$, $p<c_2'$ in (26). Then

\begin{eqnarray}
b_{c_1^{r_1} p} = Vb_{c_1^{r_1} \{ c_1 p\}_{c_2'}}, \quad 1\le p<c_2'.
\eqnum{27}
\end{eqnarray}

According to Statement 3 for the operation $Tp=\{ c_1p\}_{c_2'}$, one finds
$1\le p_0<c_2'$ and $k_0=k(p_0)>1$ such that $T^{k_0}p_0=p_0$.
Iteration of (27) then yields

\begin{eqnarray}
b_{c_1^{r_1} p_0} = V^{k_0}b_{c_1^{r_1} p_0}.
\nonumber
\end{eqnarray}
One has $b_\alpha >0$. Hence $V=1$, which is the desired result.

One is now entitled to assume that $c_2=c_1^{r_1}c_2'$ and $r_1>0$.
In that case the equations $X=Y$ are equivalent to

\begin{eqnarray}
b_p = \left\{\begin{array}{llll}
b_{\{ p\}_{c_1}}, \,\, \,\qquad \quad \qquad \qquad \qquad \qquad \mbox {\rm(a)} \\
V^{-1}b_{p/c_1}, \qquad \,\qquad \qquad \qquad \qquad \mbox {\rm (b)} \\
b_{\{ p\}_{c_2/c_1}}, \,\,\,\, \qquad \qquad \qquad \qquad \qquad \mbox {\rm (c)} \\
Ub_{p/(c_2/c_1)}, \,\, \,\, \quad \qquad \qquad \qquad \qquad \mbox {\rm (d)} \\
\end{array}\right.
\eqnum{28}
\end{eqnarray}
where $1\le p<c_2$.

This can be seen as follows. Equation (28a) can be derived
by comparing (24a) and (25a); (28b) results from a comparison
between (24b) and (25b); and (28d) from that between (24b) and (25b).
A few words are required to explain (28c). From (24b) and (25a) one has

\begin{eqnarray}
V^{-1}b_\beta = b_{\{ \beta c_1\}_{c_2}} =
b_{c_1\{ \beta \}_{(c_2/c_1)}} = V^{-1}b_{\{ \beta \}_{(c_2/c_1)}}.
\nonumber
\end{eqnarray}
The last equality in the above sequence follows from the first.

{\it Consider the case} $(c_2',c_1)=c_2'$.

The right-hand side can not be equal to $c_1$, otherwise $c_2'$ would
be divisible by $c_1$, and the number $r_1$ in the representation
$c_2=c_1^{r_1}c_2'$ will not be the maximum multiplicity of
$c_1$ in $c_2$.

To sum up, one has $c_2=c_1^{r_1}c_2'$, $r_1\ge 1$
and $c_1=(c_2')^{s_1}c_1'$, $s_1\ge 1$ where $s_1$ is the
maximum multiplicity of $c_2'$ in $c_1$. Note that $c_1'>1$, otherwise
$\ln c_1/\ln c_2$ would be rational. From (28) one has

\begin{eqnarray}
Ub_\alpha \stackrel{\rm (28d)}{=} b_{\alpha c_1^{r_1-1}c_2'}
\stackrel{\rm (28b)}{=} V^{1-r_1}b_{\alpha c_2'}
\stackrel{\rm (28a)}{=} V^{1-r_1}b_{\{ \alpha c_2'\}_{c_1}}.
\eqnum{29}
\end{eqnarray}
If $r_1=1$, then $c_2'=c_2/c_1$. Therefore,

\begin{eqnarray}
b_\alpha \stackrel{\rm (28b)}{=} V b_{\alpha c_1} =
V b_{\alpha (c_2')^{s_1}c_1'}
\stackrel{\rm (28c)}{=} V b_{\alpha c_1'}
\stackrel{\rm (28d)}{=} V b_{\{ \alpha c_1'\}_{c_2'}}.
\eqnum{30}
\end{eqnarray}
Let $(c_1',c_2')<\min (c_1',c_2')$. In that case (30) will yield $V=1$,
similarly to the above argument.

It can now be assumed that $r_1>1$. Put $\alpha =(c_2')^{s_1}\alpha'$
in (29).  One then has from (29):

\begin{eqnarray}
U b_{(c_2')^{s_1}\alpha'} =
V^{1-r_1} b_{(c_2')^{s_1}\{ c_2'\alpha'\}_{c_1'}}.
\nonumber
\end{eqnarray}
When $(c_2',c_1')<\min (c_2',c_1')$, the standard procedure would
yield $U=V^{1-r_1}$ or $V_b=V_a^{r_1}$.
However, according to Lemma 2,

\begin{eqnarray}
V_a^{\ln c_2/\ln c_1} = V_b.
\eqnum{31}
\end{eqnarray}
When $V_a\ne 1$, one has $\ln c_2/\ln c_1=r_1$, which is impossible.
Therefore, one has $V_a=1$ for the case $r_1\ge 1$ as well,
as was to be proved.

We have arrived at the situation in which $c_2=c_1^{r_1}c_2'$, \,$r\ge 1$,
\, $c_1=(c_2')^{s_1}c_1'$ and $c_1'$
is divisible by $c_2'$. We will show that a set of equations can
be written down that is similar to (28), but where $c_1$ and $c_2$
have been replaced with $c_1'$ and $c_2'$. Also, the coefficients of
$U$, $V$ which have the form $V_a^{k_1}V_b^{k_2}$ with $k_i$ integer will
be replaced with similar coefficients. Consequently, the proof will
reduce to the preceding with smaller $c_i$. The reduction process
is finite. It will terminate, when $c_1'$ is no longer divisible
by $c_2'$. Otherwise $\log c_2/\log c_1$ will be a rational number.

It remains to find the analog of (28). One has

\begin{eqnarray}
Ub_q \stackrel{\rm (28d)}{=} b_{c_2q/ c_1} =
b_{c_1^{r_1-1}c_2'q}
\stackrel{\rm (28b)}{=} (V^{-1})^{r_1-1} b_{c_2'q}.
\nonumber
\end{eqnarray}
Hence

\begin{eqnarray}
b_p = \tilde{U}b_{p/c_2'}, \quad  p<c_1c_2',
\eqnum{32}
\end{eqnarray}
where $\tilde{U}=UV^{r_1-1}$. Further,

\begin{eqnarray}
b_q \stackrel{\rm (28b)}{=} Vb_{c_1q} = V b_{(c_2')^{s_1}c_1'q}
\stackrel{\rm (32)}{=} V \tilde{U}^{s_1} b_{c_1'q},
\nonumber
\end{eqnarray}
which gives the analog of (28b):

\begin{eqnarray}
b_p = V_1^{-1}b_{p/c_1'}, \quad  p<c_1'c_2',
\eqnum{33}
\end{eqnarray}
where $V_1=V\tilde{U}^{s_1}$.

>From (32, 33) one derives the analog of (28d) with a new
coefficient  $U_1=\tilde{U}V_1$:

\begin{eqnarray}
\tilde{U} b_{p/(c_2'/c_1')}
\stackrel{\rm (32)}{=} b_{c_1'p}
\stackrel{\rm (33)}{=} V_1^{-1} b_p, \quad p<c_2'.
\eqnum{34}
\end{eqnarray}
We now are going to derive the analogs of (28a) and (28c).
>From (28a) one has

\begin{eqnarray}
b_{(c_2')^{s_1}p'} =
b_{(c_2')^{s_1}\{ p'\}_{c_1'}}, \quad p'<c_2'c_1'.
\nonumber
\end{eqnarray}
The use of (32) on both sides of the above equality will yield

\begin{eqnarray}
b_{p'} =
b_{\{ p'\}_{c_1'}}, \quad p'<c_2'c_1'.
\eqnum{35}
\end{eqnarray}
Similarly, (28c) yields

\begin{eqnarray}
b_p = b_{\{ p\}_{c_1^{r_1-1}c_2'}}.
\nonumber
\end{eqnarray}
Substituting $p=c_1^{r_1-1}p'<c_2$ and using (28b), one gets

\begin{eqnarray}
b_{p'} = b_{\{ p'\}_{c_2'}},\quad p'<c_1c_2'.
\nonumber
\end{eqnarray}
Hence one gets for $p'=c_1'q$, $q<c_2'$ with the help of (33):

\begin{eqnarray}
b_q = b_{\{ q\}_{(c_2'/c_1')}}.
\eqnum{36}
\end{eqnarray}
Relations (33-36) are the analog of (28). The proof of the theorem
is complete.

{\bf Acknowledgements.}
This work was supported by the James S. McDonnell Foundation within
the framework of the 21st Century Collaborative Activity Award for
Studying Complex Systems (project "Understanding and Prediction of
Critical Transitions in Complex Systems") and in part by the Russian
Foundation for Basic Research (Grant 02-01-00158).

\newpage

\bigskip

\centerline{\bf REFERENCES}

\begin{description}

\item
 1.\,\, Frisch U. {\it Turbulence: the Legacy of A.N. Kolmogorov}
(Cambridge University Press, 1995).

\item
 2.\,\, Kahane J.P. and Peyriere J. Sur certaines martingales
de B. Mandelbrot. {\it Adv. Math}. {\bf 22}: 131-145 (1976).

\item
 3.\,\, Holley R. and Ligget Th. Generalized potlatch and
smoothing processes. {\it Z. Wahr. Verw. Geb}. {\bf 55}: 165-195 (1981).

\item
 4.\,\, Mandelbrot B. Multiplications aleatoires et distributions
invariantes par moyenne ponderee aleatoire. {\it C.R. Acad. Sci. Paris,
Ser. A}, {\bf 278}: 289-292 and 355-358 (1974).

\item
 5.\,\, Molchan G.M. Scaling exponents and multifractal
dimensions for independent random cascades.
{\it Commun. Math. Phys}. {\bf 179}: 681-702 (1996).

\item
 6.\,\, Ossiander M. and Waymire E. Statistical estimation
for multiplicative cascades. {\it Ann. Stat}. {\bf 28}, no. 6 (2000).

\item
 7.\,\, Novikov E.A. Infinitely divisible distributions in turbulence.
{\it Phys. Rev. E} {\bf 50}: R3303-R3305 (1994).

\item
 8.\,\, Pedrizzeti G., Novikov E.A., and Praskovsky A.A. Self-similarity
and probability distributions of turbulent intermittency.
{\it Phys. Rev. E} {\bf 53}: 475-484 (1996).

\item
 9.\,\, She Z.S. and Waymire E. Quantized energy cascade and
log-Poisson statistics in fully developed turbulence.
{\it Phys. Rev. Lett}. {\bf 74}: 262-265 (1995).

\item
10.\,\, Schertzer D., Lovejoy S., Schmitt F., Chigirinskaya Y.,
and Marsan D. Multifractal cascade dynamics and turbulent
intermittency. {\it Fractals} {\bf 5}, no. 3, 427-471 (1997).

\item
11.\,\, Molchan G.M. Turbulent cascades: limitations and
a statistical test of the lognormal hypothesis. {\it Phys.
Fluids} {\bf 9}, no. 8, 2387-2396 (1997).

\item
12.\,\, Schertzer D. and Lovejoy S. Universal multifractals do exist:
comments on "A statistical analysis of mesoscale rainfall as
a random cascade". {\it J. Appl. Meteorology} {\bf 36}: 1296-1303 (1997).

\item
13.\,\, Sreenivasan K.R. and Stolovitzky G. Turbulent cascades.
{\it J. Stat. Phys}. {\bf 78}, no. 1/2, 311-333 (1995).

\item
14.\,\, Molchan G. Mandelbrot cascade measures independent of branching
parameters. {\it J. Stat. Phys}. {\bf 107}, no. 5/6, 977-988 (2002).

\item
15.\,\, Durrett R. and Ligget Th. Fixed points of smoothing
transformation. {\it Z. Wahr. Verw. Geb}. {\bf 64}: 275-301 (1983).

\end{description}

\end{document}